\documentclass{article}
\usepackage{amsfonts}
\usepackage{amsmath}
\usepackage[a4paper]{geometry}
\usepackage{amsthm}
\usepackage{amssymb}
\usepackage{enumerate}
\usepackage{mathrsfs}
\usepackage{mathtools}
\usepackage{appendix}
\usepackage{comment}
\usepackage {tikz}
\usetikzlibrary {positioning}
\usetikzlibrary{fit,backgrounds} 
\definecolor {processblue}{cmyk}{0.96,0,0,0}

\newcommand{\rn}{\mathbb{R}}
\newcommand{\fn}{\mathbb{F}}
\newcommand{\zn}{\mathbb{Z}}
\newcommand{\qn}{\mathbb{Q}}
\newcommand{\arr}{\mathcal{A}}
\newcommand{\be}{\begin{enumerate}}
\newcommand{\ee}{\end{enumerate}}
\newcommand{\bc}{\begin{center}}
\newcommand{\ec}{\end{center}}

\newcommand{\ti}{\textit}
\newcommand{\bi}{\begin{itemize}}
\newcommand{\ei}{\end{itemize}}

\theoremstyle{definition}
\newtheorem{theorem}{Theorem}
\newtheorem{lemma}[theorem]{Lemma}

\newtheorem{definition}[theorem]{Definition}

\title{Generalized threshold arrangements}
\date{}

\author{A.R.Balasubramanian  \\ \small{Chennai Mathematical Institute}\\ \small{Chennai, India}}

\begin{document}
\maketitle{}

\begin{abstract}
	An arrangement of hyperplanes is a finite collection of hyperplanes in a real Euclidean space. To such a collection one associates the characteristic polynomial that encodes the combinatorics of intersections of the hyperplanes. Finding the characteristic polynomial of the Shi threshold and the Catalan threshold arrangements was an open problem in Stanley's list of problems in \cite{stanley}. Seunghyun Seo solved both the problems by clever arguments using the finite field method in \cite{shi threshold, catalan threshold}. However, in his paper, he left open the problem of computing the characteristic polynomial of a broader class of threshold arrangements, the so-called ``generalized threshold" arrangements whose defining set of hyperplanes is given by $x_i + x_j = -l,-l+1,...,m-1,m$ for $1 \le i < j \le n$ where $l,m \in \mathbb{N}$.  In this paper, we present a method for computing the characteristic polynomial of this family of arrangements. 

\end{abstract}

\section{Introduction}

A \emph{set of hyperplane arrangements} over a vector space $V$ is a finite set of affine hyperplanes in $V$. In the literature, usually the underlying vector space is taken to be $\mathbb{R}^n$ for some $n$. To every such hyperplane arrangement, one can associate a set of \emph{regions} which are simply the connected components of the complements of the hyperplanes in the given arrangement. Over the course of their study, special interest has been given to computing the number of regions of specific families of arrangements.

To every hyperplane arrangement one associates a polynomial called the \emph{characteristic polynomial} of the arrangement. This polynomial captures lots of useful information about the underlying arrangement. In his groundbreaking work \cite{regions}, Zaslavsky proved that the number of regions of a hyperplane arrangement can be obtained by evaluating the characteristic polynomial of the arrangement at $-1$ and computing its absolute value. Hence one could calculate the number of regions provided the characteristic polynomial of the arrangement can be computed easily.

However computing the characteristic polynomial by itself is not an easy task. In his seminal paper \cite{finite field}, Athanasiadis came up with a very elegant yet powerful method for computing the characteristic polynomial of any arrangement over the rationals, the so called \emph{finite field} method. Ever since, the computation of characteristic polynomial for various arrangements have been attained through the finite field method. (See for example, \cite{catalan threshold, shi threshold}).

In the present article, we compute the characteristic polynomial of the generalized threshold arrangement which is given by the following set of hyperplanes over $\mathbb{R}^n$:

$$x_i + x_j = -l,-l+1, 1, 2, \dots, k \qquad \; 1 \le i < j \le n$$

Here $l,k$ are arbitrary non-negative integers. The case of $l = 0$ and $k = 1$, also called the Shi threshold arrangement was an open problem mentioned by Stanley in \cite{stanley}. 
In his paper, \cite{shi threshold} Seunghyun Seo had solved the characteristic polynomial for this arrangement. He also goes on to solve the case of $l=1,k=1$, the so called Catalan threshold arrangement in \cite{catalan threshold}. Seo left open the problem of computing the characteristic polynomial for the above mentioned generalized threshold arrangements.

The main contributions of this article are Theorems \ref{shithr} and \ref{catthr} which present a method for computing the characteristic polynomial of the generalized threshold arrangement. The result is obtained by an extension of the method used to compute the characteristic polynomial for the Shi and Catalan threshold arrangements by Seo.\\

\textbf{Acknowledgements: } I am extremely grateful to Prof. Priyavrat Deshpande for introducing me to the subject of hyperplane arrangements and for his invaluable suggestions on early drafts of this paper. I would also like to thank Deeparaj Bhat for useful discussions which lead to a substantial improvement in the quality of the proofs.

\section{Preliminaries}

In this section, we recall some of the basic concepts of hyperplane arrangements. The interested reader is referred to \cite{book} for more information.

\subsection{Characteristic polynomials and the finite field method}

Let $n$ be a positive integer and $K$ be a field. A \ti{hyperplane arrangement} (or an arrangement) over $K^n$ is a finite set of hyperplanes $\arr$ in $K^n$, with each hyperplane having coefficients from $K$.
\begin{definition}
	Given an arrangement $\arr$, we let $L(\arr)$ be the set of all possible intersections of hyperplanes in $\arr$, including $K^n$ and the empty set. For any two intersections $x,y$ define $x \le y$ iff $x \supseteq y$. $L(\arr)$ is called the \ti{intersection poset} of $\arr$. 	
\end{definition}

Clearly in this intersection poset, $K^n \le x$ for any element $x$. 
For any two elements $x,y \in L(\arr)$, let the interval $[x,y]$ be equal to the set $\{z : x \le z \le y\}$. Let Int$(L(\arr))$ denote the set of all intervals of the intersection poset of $\arr$. 
For the intersection poset, one can define the \ti{Mobius function} $\mu : L(\arr) \to \mathbb{Z}$ as follows:

$$\mu([x,x]) = 1 \quad \text{ and } \quad \mu([x,y]) = - \sum_{x \le z < y} \mu(y)$$ 

For any interval of the form $[K^n,x]$ we abuse notation and denote $\mu([K^n,x])$ by $\mu(x)$. Using the Mobius function, we can define the \ti{characteristic polynomial} of an arrangement $\arr$. 

\begin{definition}
	The characteristic polynomial $\chi_{\arr}(t)$ of an arrangement $\arr$ is the polynomial given by:
	$$\chi_{\arr}(t) := \sum_{x \in L(\arr)} \mu(x) t^{\dim(x)}$$
\end{definition}

Here $\dim(x)$ refers to the dimension of the affine subspace $x$ in the vector space $K^n$.

In the literature almost all hyperplanes are defined and studied over the field $K = \rn^n$. We let $\arr$ be a hyperplane arrangement over $\rn$. A \ti{region} of $\arr$ is a connected component of the complement of the hyperplanes (i.e), any connected component of $\rn^n - \bigcup_{H \in \arr} H$. Let $r(\arr)$ denote the number of regions of $\arr$.

The characteristic polynomial contains important information about the corresponding hyperplane arrangement, as evidenced by the following result. The celebrated theorem of Zaslavsky in \cite{regions} states that,

\begin{theorem}
	For an arrangement $\arr$, we have 
	$$r(\arr) = (-1)^n \chi_{\arr}(-1)$$
\end{theorem}

However the task of computing characteristic polynomials is a hard task by itself. But if $\arr$ is a hyperplane arrangement over the vector space $\qn^n$ then there is a powerful method called the \emph{finite field method} for computing the corresponding characteristic polynomial. 

Given an arrangement $\arr$ over $\qn^n$, we can get an arrangement over $\fn^n_p$ for any prime $q$ as follows: Suppose $H$ is a hyperplane in $\arr$ of the form:
$$a_1x_1 + a_2x_2 + \dots + a_nx_n = b$$

for some $a_1,a_2,\dots,a_n,b \in \qn$. Let $d$ denote the least common multiple of the denominators of $a_1,a_2,\dots,a_n,b$. Multiplying the hyperplane with $d$ on both sides, we get 
$$a_1'x_1 + a_2'x_2 + \dots + a_n'x_n = b'$$

where each $a_i' = d \cdot a_i$ and $b' = d \cdot b$. Clearly $a_1',\dots,a_n',b' \in \zn$. Let $a_i'' = a_i' \bmod p$ and $b'' = b' \bmod p$ and let $H_p$ be the hyperplane:
$$a_1''x_1 + a_2''x_2 + \dots + a_n''x_n = b''$$

For each $H \in \arr$, we get one such hyperplane $H_p$. Let the collection of all such $H_p$ be the arrangement $\arr_p$ over the vector space $\fn^n_p$. The following two theorems are well known:

\begin{theorem}
	For all but finitely many primes $p$, $L(\arr) \cong L(\arr_p)$
\end{theorem}

\begin{theorem} \cite{finite field, finite field method}
	If $L(\arr) \cong L(\arr_p)$ for some prime $p$, then
	$$\chi_{\arr}(p) = p^n - \left| \bigcup_{H \in \arr_q} H \right|$$
\end{theorem}

The above theorem is called the \ti{finite field method}. As a demonstration we use the finite field method to compute the characteristic polynomial of the \emph{braid arrangement} $\mathcal{B}_n$ given by
$$x_i - x_j = 0, \quad 1 \le i < j \le n$$

For a large enough prime $p \gg n$, we count the number of elements 
$(a_1,a_2,\dots,a_n) \in \fn^n_p$ such that $a_i - a_j \neq 0$ for any $i < j$. The number of such tuples can be clearly seen to be $p(p-1)\dots(p-n+1)$. By the finite field method we have that 
$\chi_{\mathcal{B}_n}(p) = p(p-1)\dots(p-n+1)$. Since this is true for infinitely many primes and since $\chi_{\mathcal{B}_n}$ is a polynomial in one variable, we can conclude that 
$$\chi_{\mathcal{B}_n}(t) = t(t-1)\dots(t-n+1)$$

Applying Zaslavasky's theorem one can now infer that $r(\mathcal{B}_n) = (-1)^n \cdot \chi_{\mathcal{B}_n}(-1) = n!$.\\

A few more examples relevant to our work are the Shi arrangement (\cite{shi}) $\mathcal{S}_n$, given by the following set of hyperplanes
$$x_i - x_j = 0,1 \qquad 1 \le i < j \le n$$
and the Catalan arrangement $\mathcal{C}_n$ whose defining hyperplanes are,
$$x_i - x_j = -1,0,1 \qquad 1 \le i < j \le n$$

The characteristic polynomials of $\mathcal{S}_n$ and $\mathcal{C}_n$ respectively are 
$$\chi_{\mathcal{S}_n}(t) = t(t-n)^{n-1}$$
$$\chi_{\mathcal{C}_n}(t) = t(t-n-1)(t-n-2)\dots(t-2n+1)$$

and hence the number of regions by Zaslavsky's theorem is 
$$r(\mathcal{S}_n) = (n+1)^{n-1}$$
$$r(\mathcal{C}_n) = n!C_n$$

where $C_n = \frac{1}{n+1}{2n \choose n}$ is the $n^{th}$ Catalan number.\\ 

The \ti{Shi threshold arrangement} $\mathcal{ST}_n$ is given by the following set of hyperplanes 
$$x_i + x_j = 0,1 \qquad 1 \le i < j \le n$$

and the \ti{Catalan threshold arrangement} $\mathcal{CT}_n$ by,
$$x_i + x_j = -1,0,1 \qquad 1 \le i < j \le n$$

Using the finite field method, in \cite{shi threshold} and $\cite{catalan threshold}$, Seo proved that the characteristic polynomial of $\mathcal{ST}_n$ and $\mathcal{CT}_n$ respectively is,
\begin{align*}
	\chi_{\mathcal{ST}_n}(t) &= \sum_{j \ge 0} (t-j-1)_j \; S(n,j) +
	2n \sum_{j \ge 0} (t-j-2)_j \; S(n-1,j) \\
	&+ n(n-1) \sum_{j \ge 0} (t-j-3)_j \; S(n-2,j)
\end{align*}
\begin{align*}
	\chi_{\mathcal{CT}_n}(t) &= n! \sum_{k=0}^n \sum_{l=0}^k \alpha_{n,k,l} \frac{((t-2k-1))_l}{l!}
\end{align*}

Here, 
\begin{itemize}
	\item $(x)_0 = 1$ and $(x)_k = x(x-1) \dots (x-k+1)$ for $k \ge 1$
	\item $((x))_0 = 1$ and $((x))_k = x(x-2)(x-4) \dots (x-2k+2)$ for $k \ge 1$.
	\item $S(a,b)$ denotes Stirling numbers of the second kind, (i.e) the number of ways to partition a set of $a$ objects into $b$ non-empty subsets.
	\item $\alpha_{0,0,0} = \alpha_{1,0,0} = 1$ and for $n \ge 2$ or $k^2+l^2 > 0$,
	$$\alpha_{n,k,l} = {k-1 \choose l-1}s_{n,k} + {k-2 \choose l-1}s_{n-1,k-1} + 2{k-1 \choose l}s_{n-1,k-1} + 2{k-2 \choose l}s_{n-2,k-2}$$
	where $s_{n,k} = \frac{k!}{n!}S(n,k)$
\end{itemize}

\section{Generalized threshold arrangements}

\begin{definition}
	The generalized threshold arrangement $\mathcal{T}_{n,k,l}$ is the hyperplane arrangement given by the following equations: For a fixed $k, l$ and $n$, it is given by:
	$$x_i + x_j = -l, -l+1, \dots ,0, 1, 2, \dots, k \qquad \; 1 \le i < j \le n$$	
\end{definition}

Whenever $l$ is fixed to be zero we will denote the resulting set of arrangements to be $\mathcal{ST}_{n,k}$ (i.e), $\mathcal{ST}_{n,k} := \mathcal{T}_{n,k,0}$ and whenever $l = -1,$ we will denote the resulting arrangement by $\mathcal{CT}_{n,k}$ (i.e), $\mathcal{CT}_{n,k} := \mathcal{T}_{n,k,1}$. In this section we show that to compute the characteristic polynomial of $\mathcal{T}_{n,k,l}$, it suffices to compute the characteristic polynomial of either $\mathcal{ST}_{n,k+l}$ or $\mathcal{CT}_{n,k+l+1}$.

For a prime $p$, let 
$\fn_p$ denote the finite field consisting of $p$ elements.
By the finite field method there exists a number $m$ such that if  $2r+1 > m$ and $2r+1$ is a prime then computing the characteristic polynomial $\chi_{\mathcal{T}_{n,k,l}}$ at $2r+1$, is equivalent to computing the cardinality of the following set: 
$$U = \{(x_1, x_2, \dots, x_n) \in \fn_{2r+1}^n: x_i + x_j \neq -l,-l+1,\dots,0, 1, 2, \dots, k \text{ for } i < j\}$$

Suppose $l$ is even. Consider the hyperplane arrangement $\mathcal{ST}_{n,k+l}$. By the finite field method there exists a number $m_0$ such that if $2r+1 > m_0$ and $2r+1$ is a prime then computing the characteristic polynomial $\chi_{\mathcal{ST}_{n,k+l}}$ at $2r+1$, is equivalent to computing the cardinality of the following set: 
$$U_0 = \{(x_1, x_2, \dots, x_n) \in \fn_{2r+1}^n: x_i + x_j \neq 0, 1, 2, \dots, k+l \text{ for } i < j\}$$

Let $2r+1 > \max(m,m_0)$. Consider the bijection $A : \fn_{2r+1}^n \to \fn_{2r+1}^n$ where 
$$A(x_1,x_2,\dots,x_n) = \left(x_1+\frac{l}{2},x_2+\frac{l}{2},\dots,x_n+\frac{l}{2}\right)$$ 
It is clear that $(x_1,x_2,\dots,x_n) \in U \iff A(x_1,\dots,x_n) \in U_0$. Therefore $|U| = |U_0|$. Hence, it follows that $\chi_{\mathcal{T}_{n,k,l}}(2r+1) = \chi_{\mathcal{ST}_{n,k+l}}(2r+1)$ for infinitely many primes $2r+1$. So we can conclude that $\chi_{\mathcal{T}_{n,k,l}} = \chi_{\mathcal{ST}_{n,k+l}}$. 

Suppose $l$ is odd. Consider the hyperplane arrangement $\mathcal{CT}_{n,k+l+1}$. By the finite field method there exists a number $m_1$ such that if $2r+1 > m_1$ and $2r+1$ is a prime then computing the characteristic polynomial $\chi_{\mathcal{CT}_{n,k+l+1}}$ at $2r+1$ is equivalent to computing the cardinality of the following set: 
$$U_1 = \{(x_1, x_2, \dots, x_n) \in \fn_{2r+1}^n: x_i + x_j \neq 1, 2, \dots, k+l+1 \text{ for } i < j\}$$

Let $2r+1 > \max(m,m_1)$. Consider the bijection $A' : \fn_{2r+1}^n \to \fn_{2r+1}^n$ where 
$$A'(x_1,x_2,\dots,x_n) = \left(x_1+\frac{l+1}{2},x_2+\frac{l+1}{2},\dots,x_n+\frac{l+1}{2}\right)$$
It is clear that $(x_1,x_2,\dots,x_n) \in U \iff A'(x_1,\dots,x_n) \in U_1$.
Hence we can conclude that $\chi_{\mathcal{T}_{n,k,l}} = \chi_{\mathcal{CT}_{n,k+l+1}}$. 

Combining the previous two observations, we get
\begin{lemma}
	$\chi_{\mathcal{T}_{n,k,l}} = \chi_{\mathcal{ST}_{n,k+l}}$ if $l$ is even and $\chi_{\mathcal{T}_{n,k,l}} = \chi_{\mathcal{CT}_{n,k+l+1}}$ otherwise.
\end{lemma}

\newpage

\section{Characteristic polynomial of $\mathcal{ST}_{n,k}$}
The arrangement $\mathcal{ST}_{n,k}$ is given by,
$$x_i + x_j = 0, 1, 2, \dots, k \qquad \; 1 \le i < j \le n$$

As discussed before, for a large enough prime $2r+1$ computing the characteristic polynomial $\chi_{\mathcal{ST}_{k,n}}$ at $2r+1$ is equivalent to computing the cardinality of the set: 
$$U = \{(x_1, x_2, \dots, x_n) \in \fn_{2r+1}^n: x_i + x_j \neq 0, 1, 2, \dots, k, \text{ for } i < j\}$$

Here $\fn_{2r+1} =  \{0,1,2,\dots,r,-1,-2,\dots,-r\}$ denotes the finite field with $2r+1$ elements. Instead of computing this value directly, we go via a graph-theoretic approach. Intuitively, we will form a graph with $2r+1$ vertices such that computation of the number of \emph{independent sets} of size atmost $n$ in this graph will easily enable us to count the cardinality of $U$. For general graphs, this computation is hard. However the graph that we construct will have some special properties which we will exploit to compute the required number.

Construct a graph $G = (V,E)$ where $V = \{0,1,2\dots,r,-1,-2\dots,-r\}$ and $(i,j) \in E$ iff $i+j \in \{0,1,2,\dots,k\}$. (Here $+$ refers to the addition operation in $\fn_{2r+1}$). Figure \ref{Graph} depicts the graph $G$.

An \emph{independent set} of the graph $G$ is a set of vertices $\{i_1,i_2,\dots,i_w\}$ such that for every $j$ and $l$, $(i_j,i_l) \notin E$. Similarly, a \emph{clique} of the graph $G$ is  a set of vertices $\{i_1,i_2,\dots,i_w\}$ such that for every $j$ and $l$, $(i_j,i_l) \in E$. The size of an independent set (resp. clique) is the number of vertices in it. Let $I_n$ denote the set of all independent sets of size atmost $n$ in graph $G$. It is easy to notice that if $(x_1,x_2,\dots,x_n) \in U$ then the set $\{x_1,x_2,\dots,x_n\} \in I_n$.

Consider the case when $k$ is even and $r$ is sufficiently large so that the sets $C_1 = \{0, 1, 2, \dots, k/2\}$ and $C_2 = \{-r, -r+1, \dots, -r+k/2-1\}$ are disjoint. It is easy to notice that both these sets are cliques in $G$. With this observation it is obvious that if $(x_1,x_2,\dots,x_n) \in U$, then the number of indices $i$ such that $x_i$ can be in $\{0, 1, 2, \dots, k/2\}$ is atmost one and similarly for $\{-r, -r+1, \dots, -r+k/2-1\}$. Hence if $x = (x_1,\dots,x_n) \in U$, $x$ can be belong to any of these four cases:

\begin{enumerate}
	\item There is no $i$ such that $x_i \in C_1$ or $x_i \in C_2$.
	\item There is a unique $i$ such that $x_i \in C_1$ and there is no $j$ such that $x_j \in C_2$.
	\item There is no $i$ such that $x_i \in C_1$ and there is exactly one $j$ such that $x_j \in C_2$.
	\item There is exactly one $i$ such that $x_i \in C_1$ and there is exactly one $x_j$ such that $x_j \in C_2$.
\end{enumerate}

In the remainder of this section we deal with the above cases.

\begin {center}
\begin{figure}[h]
\begin {tikzpicture}[auto ,node distance = 1.5 cm and 1.35 cm ,on grid , semithick ,
state/.style ={ circle ,top color =white , bottom color = white,
	draw, text= black , minimum width = 0.01 cm}]
\node[state, label = $0$] (0) {};
\node[state, label = $1$, node distance = 1cm] (1) [right = of 0] {};
\node[state, label = $2$] (2) [right = of 1] {};
\node[state, draw = none] (dots) [right = of 2] {$\dots \dots \dots$};

\node[state, label = $\frac{k}{2}$] (k) [right = of dots] {};
\node[state, label = $\frac{k}{2} + 1$] (k') [right = of k] {};
\node[state, draw = none, fill = none] (mdots) [right = of k'] {$\dots \dots \dots$};

\draw[dashed] (0) ++(-0.5,1) rectangle (5.5,-0.5);

\node[state, label = {\tiny $r-\frac{k}{2}+1$}] (r+k) [right = of mdots] {};
\node[state, label = {\tiny $r-\frac{k}{2}+2$}] (r+k') [right = of r+k] {};
\node[state, draw = none, fill = none] (mmdots) [right = of r+k'] {$\dots \dots \dots$};

\node[state, label = ${\scriptsize r-1}$] (r') [right = of mmdots] {};
\node[state, label = $r$] (r) [right = of r'] {};

\node[state, label = below:$-1$] (-1) [below = of 1] {};
\node[state, label = below:$-2$] (-2) [below = of 2] {};
\node[state, draw = none, fill = none] (bdots) [right = of -2] {$\dots \dots \dots$};

\node[state, label = below:$-\frac{k}{2}$] (-k) [right = of bdots] {};
\node[state, label = below:$-\frac{k}{2}-1$] (-k') [right = of -k] {};
\node[state, draw = none, fill = none] (mbdots) [right = of -k'] {$\dots \dots \dots$};

\node[state, label = below:{\tiny $-r + \frac{k}{2} - 1$}] (r-k) [right = of mbdots] {};
\node[state, label = below:{\tiny $-r + \frac{k}{2} - 2$}] (r-k') [right = of r-k] {};
\node[state, draw = none, fill = none] (mmbdots) [right = of r-k'] {$\dots \dots \dots$};

\node[state, label = below:${\scriptsize -r+1}$] (-r') [right = of mmbdots] {};
\node[state, label = below:$-r$] (-r) [right = of -r'] {};

\draw[dashed] (r-k) ++(-0.75,0.5) rectangle (15,-2.5);

\path(0) edge (0);
\path(0) edge (1);
\path(0) edge [bend left = 50] (2);
\path(0) edge [bend left = 50] (k);
\path(0) edge [bend left = 50] (k');
\path(1) edge (2);
\path(1) edge [bend left = 50] (k);
\path(2) edge [bend left = 50] (k);

\path(r-k) edge (r-k');
\path(r-k) edge [bend right = 50] (-r');
\path(r-k) edge [bend right = 50] (-r);
\path(r-k') edge [bend right = 50] (-r');
\path(r-k') edge [bend right = 50] (-r);
\path(-r') edge (-r);

\path(1) edge (-1);
\path(2) edge (-2);
\path(k) edge (-k);
\path(k') edge (-k');
\path(r+k) edge (r-k);
\path(r+k') edge (r-k');
\path(r') edge (-r');
\path(r) edge (-r);

\path(-1) edge (2);
\path(-1) edge (k);
\path(-1) edge (k');
\path(-2) edge (k);
\path(-2) edge (k');
\path(-k) edge (k');

\path(r-k) edge (r+k');
\path(r-k) edge (r');
\path(r-k) edge (r);
\path(r-k') edge (r');
\path(r-k') edge (r);
\path(-r') edge (r);

\end{tikzpicture}
\caption{Graph $G$ with the two cliques $C_1$ and $C_2$. There is an edge between $i$ and $j$ iff $i+j \in \{0,1,\dots,k\}$. The self-loops of each element in $C_1$ and $C_2$ are not present here.}
\label{Graph}
\end{figure}
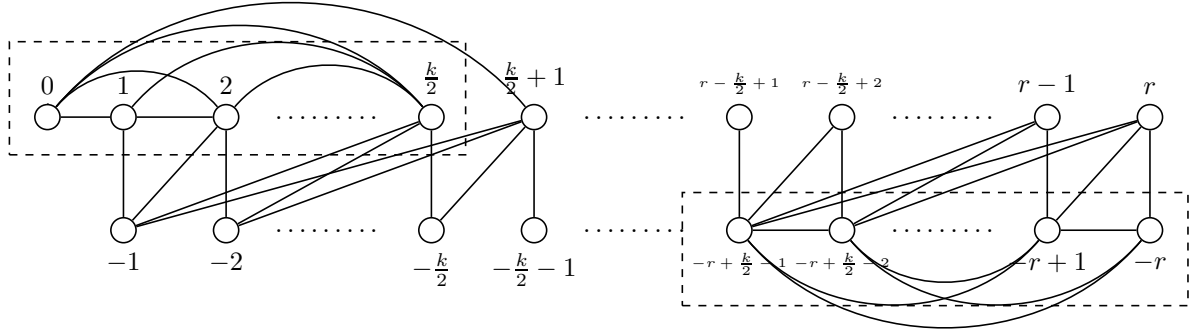
\end{center}

\subsection{No elements are picked from $C_1$ or $C_2$} 
\label{exactzero}

In this case, we would like to count the number of elements in $U$ such that they have no values from either of the cliques, (i.e) we want to count the cardinality of the set $U'$ where
$$U' = \{(x_1,x_2,\dots,x_n) : (x_1,x_2,\dots,x_n) \in U, \forall i \in [n], x_i \notin C_1 \text{ and } x_i \notin C_2\}$$

Let $G'$ be the induced subgraph of $G$ on the vertices $V \setminus (C_1 \cup C_2)$. We now show a fundamental relation between the number of independent sets in $G'$ to the number of elements in $U'$.

Suppose $I = \{y_1,y_2,\dots,y_q\}$ is a $q$-sized independent set in $G'$. We claim that $I$ contributes to $S(n,q) \cdot q!$ number of elements to the set $U'$ where $S(n,q)$ denotes Stirling numbers of the second kind. 

More precisely, we want to count the number of elements in $U'$ such that if $(x_1,x_2,\dots,x_n) \in U'$, then each $x_i = y_j$ for some $j$ \emph{and} each $y_j = x_l$ for some $l$. 
But this is precisely the number of surjective functions from $[n]$ to $[q]$ which is equal to $S(n,q) \cdot q!$.

Let the number of $q$ sized independent sets in $G'$ be denoted by $I(q)$. By the above observation, the number of elements in $U'$ is 

\begin{equation}
	N := |U'| = \sum_{1 \le q \le n} I(q) \cdot S(n,q) \cdot q! \label{first}
\end{equation}

Now for a given $q$, we derive a closed form formula for $I(q)$.

Clearly $G'$ is a bipartite graph with partitions $\{k/2+1, \dots, r\}$ and $\{-1, \dots, -r+k/2\}$. Suppose we pick $l$ elements from the first partition, say, $i_1 < i_2 < \dots < i_l$. Now we have to pick $q-l$ elements from the second partition, say, $j_1 <  j_2 < \dots < j_{q-l}$ such that the pair $(i_a,j_b)$ is not an edge for any $i_a, j_b$. In such a case, we will call the sets $\{i_1,i_2,\dots,i_l\}$ and $\{j_1,\dots,j_{q-l}\}$ \emph{compatible}. For a fixed $i_1, i_2, \dots, i_l$, let 
$f(i_1, i_2, \dots, i_l)$ be the number of such compatible sets from the second partition. Hence,

\begin{equation}
	I(q) = \sum_{0 \le l \le q} \;\; \sum_{\frac{k}{2}+1 \le i_1 < i_2 < \dots < i_l \le r} f(i_1, i_2, \dots, i_l)
	\label{indp}
\end{equation}

For a fixed $(i_1, i_2, \dots, i_l)$ with $i_1 < i_2 < \dots < i_l$, we will first compute the number of elements in the second partition such that none of the elements share an edge with any of the $i_j$. To calculate this number notice that 

\bi
\item If $a \ge -i_j$ and $(i_j,a) \notin E$, then $(i_k,a) \notin E$ for any $k \ge j$ and,
\item If $a < -i_j$, then $(i_j,a) \notin E$ for any $a$.
\ei

Hence we can cut the partition $\{-1,\dots,-r+k/2\}$ into $l+1$ slices 
\begin{multline}
	\{-1,\dots,\max \; (-r+k/2,-i_1)\}, \{-i_1-1,\dots,\max \; (-r+k/2,-i_2)\},\dots,\\
	\{-i_{l-1}-1,\dots,\max \; (-i_l,-r+k/2)\}, \{-i_l-1,\dots,-r+k/2\}	
\end{multline}
and for each $i_j$ compute the number of elements in the $j^{th}$ slice which do not share an edge with $i_j$, whose overall value will give us the number of elements which don't share an edge with any of the $i_j$. (If any of the $i_j$ is such that $-i_j \le -r+k/2$, then the set $\{-i_j-1,\dots, \max \; (-r+k/2,-i_{j+1})\}$ is taken to be empty).

More formally

\begin{definition}
	Suppose $x$ belongs to the $j^{th}$ slice of $\{-1,\dots,-r+k/2\}$, (i.e), $x \in \{-i_{j-1}-1,\dots,\max \; (-i_j,-r+k/2)\}$. Then $x$ is called \emph{independent} if $x$ does not share an edge with $i_j$.
\end{definition} 
We have the following two cases:

\bi
\item Suppose $i_l < r-k/2$. In this case, $-i_l-1 > -r+k/2$. 
The number of elements in the first slice which are independent is clearly $\max \; (i_1 - (k+1),0)$. It is also easy to see that the number of elements in the $j^{th}$ slice which are independent is $\max \; (i_j - i_{j-1} - (k+1),0)$. Finally all the elements in $\{-i_l-1,\dots,-r+k/2\}$ don't share an edge with any of the $i_j$'s and so the total number of elements in the second partition of $G'$ which don't share an edge with any of the $i_j$'s is 
\begin{multline*}
	\max \; (i_1 - (k+1),0) + \max \; (i_2 - i_1 - (k+1),0) + \dots \\+ \max \; (i_l-i_{l-1}- (k+1),0) + \max \; (r-k/2-i_l,0)	
\end{multline*}

\item Suppose $\exists \; l'$ s.t. $i_l' \ge r-k/2$. In this case notice that for any $l'' \ge l'$, $i_{l''} \ge r-k/2$. Similar to the first case, it is clear that for any $j < l'$, the number of independent elements in the $j^{th}$ slice is $\max \; (i_j - i_{j-1} - (k+1),0)$. Notice that for any $j > l'$, the $j^{th}$ slice is empty and so they don't contribute any elements. All it suffices now is to give the number of independent elements in the $l'^{th}$ slice. This is easily seen to be $\{\max(i_{l-1} - i_l - (k+1),0)\}$. Noticing that for any $j > l'$, the expression $\max(i_j - i_{j-1} - (k+1),0)$ always gives $0$ as the answer, the number of independent elements in this case as well is
\begin{multline*}
	\max \; (i_1 - (k+1),0) + \max \; (i_2 - i_1 - (k+1),0) + \dots \
	\\ + \max \; (i_l-i_{l-1}- (k+1),0) + \max \; (r-k/2-i_l,0)	
\end{multline*}
\ei

Hence the total number of independent elements is 
\begin{multline}
\max \; (i_1 - (k+1),0) + \max \; (i_2 - i_1 - (k+1),0) + \dots \
\\ + \max \; (i_l-i_{l-1}- (k+1),0) + \max \; (r-k/2-i_l,0)	
\end{multline}
We will denote this value by $g(i_1, \dots, i_l)$. It is then easily seen that $f(i_1, \dots, i_l) = {g(i_1, \dots, i_l) \choose q-l}$. 

\subsection{Exactly one element is picked from either $C_1$ or $C_2$}
\label{exactone}

We consider two different cases.

Let $0 \le p \le k/2$ be the chosen element from the first clique. In this case, we would like to count the number of elements in $U$ s.t. it contains $p$ in some position and contains no elements from the second clique. Let $U'_p \subseteq U$ be such that all elements in $U'_p$ contain $p$ in some position and contain no elements from the second clique, (i.e),
$$U'_p = \{(x_1,\dots,x_n) : (x_1,\dots,x_n) \in U, \; \exists i \text { such that } x_i = p, \text{ and } \forall j \; x_j \notin C_2\}$$

Let $F = \{0, 1, 2, \dots, k-p\} \cup \{-1,-2,\dots,-p\} \cup \{-r, -r+1, \dots, -r+k/2-1\}$ and let $G'_p$ be the induced subgraph of $G$ on the vertex set $V \setminus F$. Once again, we can relate the independent sets in $G'_p$ to the elements in $U'_p$ as follows:

Suppose $I = \{y_1,\dots,y_q\}$ is a $q$-sized independent set in $G'_p$. By a similar argument as above, it is clear that $I$ contributes to $S(n-1,q) \cdot q!$ number of elements to the set $U'_p$. Let the number of $q$ sized independent sets in $G'_p$ be denoted by $I_p(q)$. Hence, the number of elements in $U'_p$ is 

\begin{equation}
	N_p := |U_p'| = n \left(\sum_{1 \le q \le n-1	} I_p(q) \cdot S(n-1,q) \cdot q! \right) \label{second}	
\end{equation}

The graph $G'_p$ is a bipartite graph with partitions $\{k-p+1,\dots,r\}$ and $\{-p-1,\dots,-r+k/2\}$.
Similar to the previous section, we can show that

$$I_p(q) = \sum_{0 \le l \le q} \;\; \sum_{k-p+1 \le i_1 < i_2 < \dots < i_l \le r} f(i_1,i_2,\dots,i_l)$$

where $f(i_1,\dots,i_l)$ is the number of compatible sets for $(i_1,i_2,\dots,i_l)$ in the second partition of the graph $G'_p$. To compute this, once again we compute the number of independent elements in each of the slices 
\begin{multline*}
	\{-p-1,\dots,-i_1,\},
	\{-i_1-1,\dots,-i_2\},\dots,
	\{-i_{l-1}-1,\dots,-i_l\},\{-i_l-1,\dots,-r+k/2\}	
\end{multline*}
which we will denote by $g(i_1,\dots,i_l)$. It can once again be shown that the number of independent elements is 
\begin{multline*}
	\max \; (i_1- p -(k+1),0) + \max \; (i_2-i_1-(k+1),0) + \max \; (i_3-i_2-(k+1),0) + \dots \\
	+ \max \; (i_l-i_{l-1}-(k+1),0) + \max \; (r-k/2-i_l,0)	
\end{multline*}

We can do a similar computation for the case when $p$ is chosen from the second clique.
\subsection{Exactly one element is picked from both $C_1$ and $C_2$}
\label{exacttwo}

Let $0 \le p_1 \le k/2$ and $0 \le p_2 \le k/2-1$ be such that $p_1, -r + p_2$ are the two chosen elements from each of the cliques. Let $U'_{p_1,p_2}$ be the subset of $U$ such that each element in $U'_{p_1,p_2}$ contains both $p_1$ and $-r + p_2$ in some positions.
Let $F = \{0, 1, 2, \dots, k-p_1\} \cup \{r-p_2, r-p_2+1, \dots, r\} \cup \{-1,-2,\dots,-p_1\} \cup \{-r, -r+1, \dots, -r+k-(p_2+1)\}$. Let $G'_{p_1,p_2}$ be the induced subgraph of $G$ on the vertices $V \setminus F$. Once again, we can relate the independent sets in $G'_{p_1,p_2}$ to the elements in $U'_{p_1,p_2}$ as before and compute $N_{p_1,p_2}$, where

\begin{equation}
	N_{p_1,p_2} := |U'_{p_1,p_2}| = n(n-1) \left(\sum_{1 \le q \le n-2} I_{p_1,p_2}(q) \cdot S(n-2,q) \cdot q! \right) \label{third}	
\end{equation}

where the function $I_{p_1,p_2}(q)$ is the number of independent sets of size $q$ in $G'_{p_1,p_2}$ which can be computed as before.\\

Combining equations \ref{first}, \ref{second} and \ref{third} we get,
\begin{multline}
	|U| = \chi_{\mathcal{ST}_{k,n}}(2r+1) = N + \sum_{0 \le p \le \frac{k}{2}} N_p + \sum_{-r \le p \le -r+\frac{k}{2}-1} N_p + \sum_{0 \le p_1 \le \frac{k}{2}, 0 \le p_2 \le \frac{k}{2}-1} N_{p_1,p_2}	
\end{multline}

Therefore, we have
\begin{theorem} \label{shithr}
	The characterisitic polynomial of $\mathcal{ST}_{k,n}$ at any sufficiently large odd number $t$ is given by,
	\begin{align*}
		\chi_{\mathcal{ST}_{k,n}}(t) &= \sum_{1 \le q \le n} I(q) \cdot S(n,q) \cdot q! 
		+\sum_{0 \le p \le \frac{k}{2}} n \left( \sum_{1 \le q \le n-1} I_p(q) \cdot S(n-1,q) \cdot q!\right)\\
		&+ \sum_{\lfloor -\frac{t}{2} \rfloor \le p \le \lfloor -\frac{t}{2} \rfloor +\frac{k}{2}-1} n \left( \sum_{1 \le q \le n-1} I_p(q) \cdot S(n-1,q) \cdot q!\right)\\
		&+ \sum_{0 \le p_1 \le \frac{k}{2}, 0 \le p_2 \le \frac{k}{2}-1} n(n-1) \left( \sum_{1 \le q \le n-2} I_{p_1,p_2}(q) \cdot S(n-2,q) \cdot q! \right)
	\end{align*}
\end{theorem}

Suppose $k$ is odd. Then notice that in the original graph $G$, we would have the following two cliques:
$\{0, 1, 2, \dots, (k-1)/2\}$ and $\{-r, -r+1, \dots, -r+(k+1)/2-1\}$. We can proceed as before and perform a similar computation to calculate the characteristic polynomial. The final expression we would get in such a case is,
\begin{align*}
\chi_{\mathcal{ST}_{k,n}}(t) &= \sum_{1 \le q \le n} I(q) \cdot S(n,q) \cdot q! 
+ \sum_{0 \le p \le \frac{k-1}{2}} n \left( \sum_{1 \le q \le n-1} I_p(q) \cdot S(n-1,q) \cdot q!\right)\\
&+ \sum_{\lfloor -\frac{t}{2} \rfloor \le p \le \lfloor -\frac{t}{2} \rfloor +\frac{k+1}{2}-1} n \left( \sum_{1 \le q \le n-1} I_p(q) \cdot S(n-1,q) \cdot q!\right)\\
&+ \sum_{0 \le p_1 \le \frac{k-1}{2}, 0 \le p_2 \le \frac{k+1}{2}-1} 
 n(n-1) \left( \sum_{1 \le q \le n-2} I_{p_1,p_2}(q) \cdot S(n-2,q) \cdot q! \right)
\end{align*} 

\section{Characteristic polynomial of $\mathcal{CT}_{n,k}$}
The arrangement $\mathcal{CT}_{n,k}$ is given by,
$$x_i + x_j = 1, 2, \dots, k \qquad \; 1 \le i < j \le n$$

For a large enough prime $2r+1$, computing the characteristic polynomial $\chi_{\mathcal{ST}_{k,n}}$ at $2r+1$ is equivalent to computing the cardinality of the set: 
$$U' = \{(x_1, x_2, \dots, x_n) \in \fn_{2r+1}^n: x_i + x_j \neq 1, 2, \dots, k, \text{ for } i < j\}$$

As before we construct a graph $G' = (V',E')$ where $V' = \fn_{2r+1}$ and $(i,j) \in E'$ iff $i+j \in \{1,2,\dots,k\}$. Consider the case when $k$ is even and $r$ is sufficiently large so that $C_1' = \{1,2,\dots,k/2\}$ and $C_2' = \{-r,-r+1,\dots,-r+k/2-1\}$ are disjoint. It is easy to see that $C_1'$ and $C_2'$ are cliques and so if $(x_1,\dots,x_n) \in U'$ there exists atmost one index $i$ such that $x_i \in C_1'$ and atmost one index $j$ such that $x_j \in C_2'$.
Further we also have that $0$ is adjacent to every element in $C_1'$, but not to itself. Therefore the computation naturally splits into six cases:

\begin{enumerate}
	\item There is no $i$ such that $x_i \in C_1'$ or $x_i \in C_2'$ or $x_i = 0$.
	\item There is no $i$ such that $x_i \in C_1'$ or $x_i \in C_2'$ and there is atleast one $j$ such that $x_j = 0$.
	\item There is a unique $i$ such that $x_i \in C_1'$ and there is no $j$ such that $x_j \in C_2'$ and there is no $k$ such that $x_k = 0$. 
	\item There is no $i$ such that $x_i \in C_1'$ and there is exactly one $j$ such that $x_j \in C_2'$ and there is no $k$ such that $x_k = 0$.
	\item There is no $i$ such that $x_i \in C_1'$ and there is exactly one $j$ such that $x_j \in C_2'$ and there is atleast one $k$ such that $x_k = 0$.
	\item There is exactly one $i$ such that $x_i \in C_1'$ and there is exactly one $x_j$ such that $x_j \in C_2'$ and there is no $k$ such that $x_k = 0$.
\end{enumerate}

\begin {center}
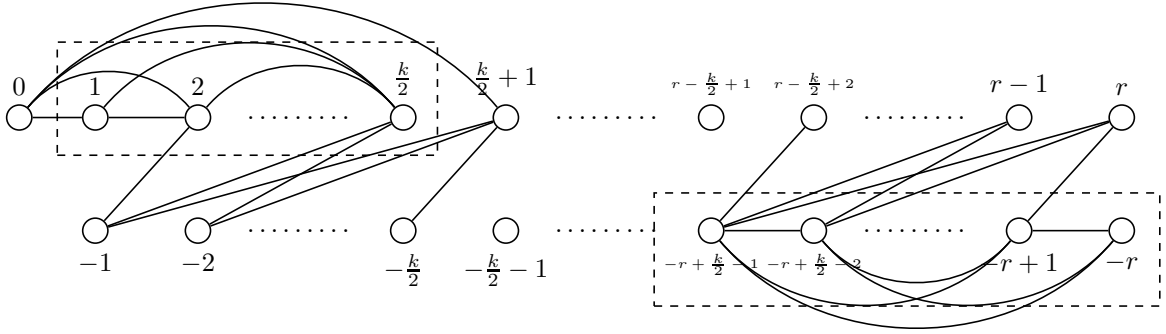
\begin{figure}[h]
	\begin {tikzpicture}[auto ,node distance = 1.5 cm and 1.35 cm ,on grid , semithick ,
	state/.style ={ circle ,top color =white , bottom color = white,
		draw, text= black , minimum width = 0.01 cm}]
	\node[state, label = $0$] (0) {};
	\node[state, label = $1$, node distance = 1cm] (1) [right = of 0] {};
	\node[state, label = $2$] (2) [right = of 1] {};
	\node[state, draw = none] (dots) [right = of 2] {$\dots \dots \dots$};

	\node[state, label = $\frac{k}{2}$] (k) [right = of dots] {};
	\node[state, label = $\frac{k}{2} + 1$] (k') [right = of k] {};
	\node[state, draw = none, fill = none] (mdots) [right = of k'] {$\dots \dots \dots$};
	
	\draw[dashed] (1,0) ++(-0.5,1) rectangle (5.5,-0.5);
	
	\node[state, label = {\tiny $r-\frac{k}{2}+1$}] (r+k) [right = of mdots] {};
	\node[state, label = {\tiny $r-\frac{k}{2}+2$}] (r+k') [right = of r+k] {};
	\node[state, draw = none, fill = none] (mmdots) [right = of r+k'] {$\dots \dots \dots$};
	
	\node[state, label = ${\scriptsize r-1}$] (r') [right = of mmdots] {};
	\node[state, label = $r$] (r) [right = of r'] {};
	
	\node[state, label = below:$-1$] (-1) [below = of 1] {};
	\node[state, label = below:$-2$] (-2) [below = of 2] {};
	\node[state, draw = none, fill = none] (bdots) [right = of -2] {$\dots \dots \dots$};
	
	\node[state, label = below:$-\frac{k}{2}$] (-k) [right = of bdots] {};
	\node[state, label = below:$-\frac{k}{2}-1$] (-k') [right = of -k] {};
	\node[state, draw = none, fill = none] (mbdots) [right = of -k'] {$\dots \dots \dots$};

	\node[state, label = below:{\tiny $-r + \frac{k}{2} - 1$}] (r-k) [right = of mbdots] {};
	\node[state, label = below:{\tiny $-r + \frac{k}{2} - 2$}] (r-k') [right = of r-k] {};
	\node[state, draw = none, fill = none] (mmbdots) [right = of r-k'] {$\dots \dots \dots$};

	\node[state, label = below:${\scriptsize -r+1}$] (-r') [right = of mmbdots] {};
	\node[state, label = below:$-r$] (-r) [right = of -r'] {};
	
	\draw[dashed] (r-k) ++(-0.75,0.5) rectangle (15,-2.5);

	\path(0) edge (1);
	\path(0) edge [bend left = 50] (2);
	\path(0) edge [bend left = 50] (k);
	\path(0) edge [bend left = 50] (k');
	\path(1) edge (2);
	\path(1) edge [bend left = 50] (k);
	\path(2) edge [bend left = 50] (k);
	
	\path(r-k) edge (r-k');
	\path(r-k) edge [bend right = 50] (-r');
	\path(r-k) edge [bend right = 50] (-r);
	\path(r-k') edge [bend right = 50] (-r');
	\path(r-k') edge [bend right = 50] (-r);
	\path(-r') edge (-r);
	
	
	\path(-1) edge (2);
	\path(-1) edge (k);
	\path(-1) edge (k');
	\path(-2) edge (k);
	\path(-2) edge (k');
	\path(-k) edge (k');
	
	\path(r-k) edge (r+k');
	\path(r-k) edge (r');
	\path(r-k) edge (r);
	\path(r-k') edge (r');
	\path(r-k') edge (r);
	\path(-r') edge (r);

\end{tikzpicture}
\caption{Graph $G$ with the two cliques $C_1'$ and $C_2'$. There is an edge between $i$ and $j$ iff $i+j \in \{1,\dots,k\}$. The self-loops of each element in $C_1'$ and $C_2'$ are not present here.}
\label{Graph1}
\end{figure}
\end{center}

\subsubsection*{The first and second cases}
We show how to compute the number of elements satisfying case 1. Let 
$$U'' = \{(x_1,x_2,\dots,x_n) : (x_1,x_2,\dots,x_n) \in U', \forall i \in [n], x_i \notin C_1' \text{ and } x_i \notin C_2' \text{ and } x_i \neq 0\}$$

Let $G''$ be the induced subgraph on $V \setminus (C_1'\cup C_2' \cup \{0\})$ and let $I'(q)$ denote the number of independent sets of size $q$ in $G''$. Similar to section \ref{exactzero}, we have
\begin{equation*}
|U''| = \sum_{1 \le q \le n} I'(q) \cdot S(n,q) \cdot q! 
\end{equation*}
If we define $f'(i_1,i_2,\dots,i_l)$ to be the number of compatible sets for $\{i_1,\dots,i_l\}$ in the second partition of $G''$, then we get 
\begin{equation}
I'(q) = \sum_{0 \le l \le q} \;\; \sum_{\frac{k}{2}+1 \le i_1 < i_2 < \dots < i_l \le r} f'(i_1, i_2, \dots, i_l)
\label{succinct}
\end{equation}

By a similar reasoning as before, we can easily conclude that,
$f'(i_1,\dots,i_l) = {g'(i_1,\dots,i_l) \choose q-l}$ where $g'(i_1,\dots,i_l)$ is 
\begin{multline*}
\max \; ((i_1-1) - (k+1),0) + \max \; ((i_2 - i_1-1) - (k+1),0) + \dots \
\\ + \max \; ((i_l-i_{l-1}-1)- (k+1),0) + \max \; (r-k/2-i_l+1,0)	
\end{multline*}

We will now show how to compute the number of elements satisfying case 2. We want to compute the cardinality of the set $U'' \subseteq U'$ given by:
$$U'' = \{(x_1,\dots,x_n) : (x_1,\dots,x_n) \in U', \forall i \; x_i \notin C_1', x_i \notin C_2', \; \exists j \; x_j = 0\}$$

To do this we fix a number $l$ such that $1 \le l \le n$ and let $h(l)$ denote the set of all elements in $U''$ such that there exists positions $i_1 < i_2 < \dots < i_l$ such that
$x_j = 0 \iff j \in \{i_1,i_2,\dots,i_l\}$. Therefore we have
$$|U''| = \sum_{1 \le l \le n} {n \choose l} h(l)$$ 

Hence it suffices to compute the number $h(l)$ for each $l$. Consider the induced subgraph $G''$ on the vertex set $V \setminus (C_1' \cup C_2' \cup \{0\})$. Letting $I'(q)$ be the number of independent sets of size $q$ in $G''$, we can easily conclude that
$$h(l) = \sum_{1 \le q \le n-l} I'(q) \cdot S(n-l,q) \cdot q!$$

By equation \ref{succinct} we already know how to compute $I'(q)$ and so we are done.

\subsubsection*{The third and fourth cases}

Consider the set of elements which satisfy the third case. Let $1 \le p \le k/2$ be the element chosen from the first clique. Define
$$U''_p = \{(x_1,\dots,x_n) : (x_1,\dots,x_n) \in U', \; \exists i \text { such that } x_i = p, \text{ and } \forall j \; x_j \notin C_2', x_j \neq 0\}$$

Let $F = \{0,1,\dots,k-p\} \cup \cup \{-1,-2,\dots,-p+1\} \cup \{-r, -r+1, \dots, -r+k/2-1\}$ and let $G''_p$ be the induced subgraph of $G'$ on $V' \setminus F$. Letting $I'_p(q)$ denote the number of independent sets of size $q$ in $G''_p$ we get that
\begin{equation*}
|U_p''| = n \left(\sum_{1 \le q \le n-1	} I'_p(q) \cdot S(n-1,q) \cdot q! \right)
\end{equation*}

Similar to section \ref{exactone}, we can write
$$I'_p(q) = \sum_{0 \le l \le q} \;\; \sum_{k-p+1 \le i_1 < i_2 < \dots < i_l \le r} f'(i_1,i_2,\dots,i_l)$$

where $f'(i_1,\dots,i_l)$ is the number of compatible sets for $(i_1,i_2,\dots,i_l)$ in the second partition of the graph $G''_p$. Letting $g'(i_1,\dots,i_l)$ equal to
\begin{multline*}
\max \; ((i_1- p -1) -(k+1),0) + \max \; ((i_2-i_1-1)-(k+1),0) + \max \; ((i_3-i_2-1)-(k+1),0) + \dots \\
+ \max \; ((i_l-i_{l-1}-1)-(k+1),0) + \max \; (r-k/2-i_l+1,0)	
\end{multline*}
we get that $f'(i_1,\dots,i_l) = {g'(i_1,\dots,i_l) \choose q-l}$

A similar computation can be done when the fourth case is satisfied (i.e) $p$ is picked from the second clique and is of the form $-r \le p \le -r+k/2-1$.

\subsubsection*{The fifth and the sixth cases}

Let $0 \le p \le k/2-1$ such that $-r+p$ is the chosen element from the second clique. Define
\begin{multline*}
	U''_{p,0} = \{(x_1,\dots,x_n) : (x_1,\dots,x_n) \in U', \; \exists i \text { such that } x_i = -r+p, \text{ and } \forall j \; x_j \notin C_1', \; \\ \text{ and } \exists k \text{ such that } x_k = 0\}	
\end{multline*}

Similar to the second case, we fix a number $l$ such that $1 \le l \le n-1$ and let $h(l)$ denote the set of all elements in $U''_{p,0}$ such that there exists positions $i_1 < i_2 < \dots < i_l$ such that
$x_j = 0 \iff j \in \{i_1,i_2,\dots,i_l\}$. Therefore we have
$$|U''_{p,0}| = n \left(\sum_{1 \le l \le n-1} {n \choose l} h(l)\right)$$ 

Let $F= \{r-p+1,\dots,r\} \cup \{-r,-r+1,\dots,-r+k-(p+1)\} \cup \{0,1,\dots,k\}$ and consider the induced subgraph $G'$ on $V' \setminus F$. It can then be seen that,
$$h(l) = \sum_{1 \le q \le n-l-1} I'_p(q) \cdot S(n-l-1,q) \cdot q!$$

where $I'_p(q)$ is the number of independent sets of size $q$ in $G''_p$, which we know how to compute.\\

For the sixth case, for every $p_1,p_2$ such that $0 \le p_1 \le k/2, 0 \le p_2 \le k/2-1$, we can once again define a function $I'_{p_1,p_2}(q)$ similar to the one defined in section \ref{exacttwo}, using which we can compute the number of required elements.\\

\begin{theorem} \label{catthr}
	Combining all the previous cases together, when $k$ is even we get,
	\begin{align*}
	\chi_{\mathcal{CT}_{k,n}}(t) &= \sum_{1 \le q \le n} I'(q) \cdot S(n,q) \cdot q! + \sum_{1 \le l \le n} {n \choose l} \sum_{1 \le q \le n-l} I'(q) \cdot S(n-l,q) \cdot q!\\
	&+\sum_{1 \le p \le \frac{k}{2}} n \left(\sum_{1 \le q \le n-1} I'_p(q) \cdot S(n-1,q) \cdot q!\right)\\
	&+ \sum_{\lfloor -\frac{t}{2} \rfloor \le p \le \lfloor -\frac{t}{2} \rfloor +\frac{k}{2}-1} n \left( \sum_{1 \le q \le n-1} I'_p(q) \cdot S(n-1,q) \cdot q!\right)\\
	&+\sum_{\lfloor -\frac{t}{2} \rfloor \le p \le \lfloor -\frac{t}{2} \rfloor +\frac{k}{2}-1} n \left( \sum_{1 \le l \le n-1} {n \choose l} \sum_{1 \le q \le n-l-1} I'_p(q) \cdot S(n-l-1,q) \cdot q! \right)\\
	&+ \sum_{0 \le p_1 \le \frac{k}{2}, 0 \le p_2 \le \frac{k}{2}-1} 
	n(n-1) \left(\sum_{1 \le q \le n-2} I'_{p_1,p_2}(q) \cdot S(n-2,q) \cdot q! \right)
	\end{align*}
	
	For odd $k$ the expression is,
	\begin{align*}
	\chi_{\mathcal{CT}_{k,n}}(t) &= \sum_{1 \le q \le n} I'(q) \cdot S(n,q) \cdot q! + \sum_{1 \le l \le n} {n \choose l} \sum_{1 \le q \le n-l} I'(q) \cdot S(n-l,q) \cdot q!\\
	&+\sum_{1 \le p \le \frac{k-1}{2}} n \left(\sum_{1 \le q \le n-1} I'_p(q) \cdot S(n-1,q) \cdot q!\right)\\
	&+ \sum_{\lfloor -\frac{t}{2} \rfloor \le p \le \lfloor -\frac{t}{2} \rfloor +\frac{k+1}{2}-1} n \left( \sum_{1 \le q \le n-1} I'_p(q) \cdot S(n-1,q) \cdot q!\right)\\
	&+\sum_{\lfloor -\frac{t}{2} \rfloor \le p \le \lfloor -\frac{t}{2} \rfloor +\frac{k+1}{2}-1} n \left( \sum_{1 \le l \le n-1} {n \choose l} \sum_{1 \le q \le n-l-1} I'_p(q) \cdot S(n-l-1,q) \cdot q! \right)\\
	&+ \sum_{0 \le p_1 \le \frac{k-1}{2}, 0 \le p_2 \le \frac{k+1}{2}-1} 
	n(n-1) \left(\sum_{1 \le q \le n-2} I'_{p_1,p_2}(q) \cdot S(n-2,q) \cdot q! \right)
	\end{align*}
	
\end{theorem}

\section{Some computations}

In this subsection, we present the characteristic polynomials for $\mathcal{ST}_{k,n}$ for some values of $k,n$. In particular, our computation agrees with the polynomials presented in \cite{shi threshold}. Further there was an error in the computation in \cite{shi threshold} for the characteristic polynomial of $\mathcal{ST}_{1,2}(t)$, which we have corrected here.

\begin{align*}
\chi_{\mathcal{ST}_{1,2}}(t) &= t^2 - 2t - 2\\
\chi_{\mathcal{ST}_{1,3}}(t) &= t^3 - 6t^2 + 12t -8\\
\chi_{\mathcal{ST}_{1,4}}(t) &= t^4 - 12t^3 + 60t^2 -142t + 130\\
\chi_{\mathcal{ST}_{1,5}}(t) &= t^5 - 20t^4 + 180t^3 - 870t^2 + 2190t -2252 \\
\chi_{\mathcal{ST}_{2,3}}(t) &= t^3 - 3t^2 - 3t -36\\
\chi_{\mathcal{ST}_{2,4}}(t) &= t^4 - 6t^3 + 27t^2 - 663t + 3583\\
\chi_{\mathcal{ST}_{2,5}}(t) &= t^5 - 10t^4 + 125t^3 - 4125t^2 + 49410t - 188578\\
\chi_{\mathcal{ST}_{3,4}}(t) &= t^4 - 204t^2 +1240t + 740\\
\chi_{\mathcal{ST}_{3,5}}(t) &= t^5 - 420t^3 + 2730t^2 + 44860t - 425199\\
\end{align*}

Applying Zaslavsky's theorem, we get that the sequence $r(\mathcal{ST}_{1,n})_{n \ge 2}$ is,

$$1,27,345,5513,\dots$$

The sequence $r(\mathcal{ST}_{2,n})_{n \ge 3}$ is,

$$37, 4280, 242249 \dots$$

\end{document}